\documentclass{article}

\usepackage{amsmath}
\usepackage{amssymb}

\newtheorem{thm}{Theorem}
\newtheorem{conj}{Conjecture}
\newtheorem{lem}{Lemma}
\newtheorem{cor}{Corollary}
\newtheorem{ques}{Question}

\title{Approximating reals by sums of rationals}
\author{Tsz Ho Chan}

\begin{document}
\maketitle
\begin{abstract}
We consider the question of approximating any real number $\alpha$
by sums of $n$ rational numbers $\frac{a_1}{q_1} + \frac{a_2}{q_2} +
... + \frac{a_n}{q_n}$ with denominators $1 \leq q_1, q_2, ... , q_n
\leq N$. This leads to an inquiry on approximating a real number by
rational numbers with a prescribed number of prime factors in the
denominator.
\end{abstract}

%------------------------------------------------------------------------
\section{Introduction and main results}
In [\ref{C}], the author generalized Dirichlet's diophantine
approximation theorem and studied the question of approximating any
real number by a sum of two rational numbers.
\begin{ques} \label{q1}
Find a good upper bound for
\begin{equation*}
\Big| \alpha - \frac{a_1}{q_1} - \frac{a_2}{q_2} \Big|
\end{equation*}
with integers $a_1, a_2$ and $1 \leq q_1, q_2 \leq N$.
\end{ques}
It turns out that one way to look at the question is to compare it
with single rational approximations of $\alpha$. From [\ref{C}], we
have
\begin{thm} \label{thm0}
For any $\epsilon > 0$ and any $N \geq 1$, suppose $\alpha$ has a
rational approximation $|\alpha - \frac{a}{q}| \leq \frac{1}{q
N^{3/2}}$ for some integers $a$, $1 \leq q \leq N^{3/2}$ and $(a,q)
= 1$. Then
\begin{equation*}
\Big|\alpha - \frac{a_1}{q_1} - \frac{a_2}{q_2}\Big| \ll_\epsilon
\frac{1}{q N^{3/2 - \epsilon}}
\end{equation*}
for some integers $a_1$, $a_2$, $1 \leq q_1, q_2 \leq N$.
\end{thm}
It was conjectured in [\ref{C}] that
\begin{conj} \label{conj1}
For any small $\epsilon > 0$ and any $N \geq 1$, suppose $\alpha$
has a rational approximation $|\alpha - \frac{a}{q}| \leq \frac{1}{q
N^2}$ for some integers $a$, $1 \leq q \leq N^2$ and $(a,q) = 1$.
Then
$$\Big|\alpha - \frac{a_1}{q_1} - \frac{a_2}{q_2}\Big| \ll_\epsilon
\frac{1}{q N^{2 - \epsilon}}$$ for some integers $a_1$, $a_2$, $1
\leq q_1, q_2 \leq N$.
\end{conj}

More generally, one can consider the following
\begin{ques} \label{q2}
Find a good upper bound for
\begin{equation*}
\Big| \alpha - \frac{a_1}{q_1} - \frac{a_2}{q_2} - ... -
\frac{a_n}{q_n} \Big|
\end{equation*}
with integers $a_1, a_2, ..., a_n$ and $1 \leq q_1, q_2, ..., q_n
\leq N$.
\end{ques}

Towards it, we have
\begin{thm} \label{thm1}
Let $\epsilon > 0$ and $N \geq 1$ be any real numbers. For any
positive integer $n \leq \frac{\epsilon \log N}{6 \log \log N}$, let
$$\kappa(n) := \frac{3n}{4} - \frac{[n/3] + 1}{4}.$$
Suppose $\alpha$ has a rational approximation $|\alpha -
\frac{a}{q}| \leq \frac{1}{q N^{\kappa(n)}}$ for some integers $a$,
$1 \leq q \leq N^{\kappa(n)}$ and $(a,q) = 1$. Then
\begin{equation} \label{1}
\Big|\alpha - \frac{a_1}{q_1} - \frac{a_2}{q_2} - ... -
\frac{a_n}{q_n} \Big| \ll_\epsilon \frac{1}{q N^{\kappa(n) -
\epsilon}}
\end{equation}
for some integers $a_1$, $a_2$, ..., $a_n$ and distinct prime
numbers $1 \leq q_1, q_2, ..., q_n \leq N$.
\end{thm}

When $n=3$, one can use an idea from [\ref{C}] to get a better
exponent.
\begin{thm} \label{thm2}
Let $\epsilon, \epsilon' > 0$ and $N \geq 1$ be any real numbers.
Suppose $\alpha$ has a rational approximation $|\alpha -
\frac{a}{q}| \leq \frac{1}{q N^{2 - \epsilon'}}$ for some integers
$a$, $1 \leq q \leq N^{2 - \epsilon'}$ and $(a,q) = 1$. Then
\begin{equation} \label{2}
\Big|\alpha - \frac{a_1}{q_1} - \frac{a_2}{q_2} - \frac{a_3}{q_3}
\Big| \ll_\epsilon \frac{1}{q N^{2 - \epsilon' - \epsilon}}
\end{equation}
for some integers $a_1$, $a_2$, $a_3$ and $1 \leq q_1, q_2, q_3 \leq
N$.
\end{thm}

%Note: Igor Shparlinski communicated to the author that the exponent
%$\kappa(n)$ can be improved to $3/2$ instead of $7/4$ (from above)
%when $n=3$ by a different method.

Similar to Conjecture \ref{conj1}, we have
\begin{conj} \label{conj2}
Let $\epsilon > 0$ and $N \geq 1$ be any real numbers. Let $n \leq
\epsilon \log N$ be a positive integer. Suppose $\alpha$ has a
rational approximation $|\alpha - \frac{a}{q}| \leq \frac{1}{q N^n}$
for some integers $a$, $1 \leq q \leq N^n$ and $(a,q) = 1$. Then
\begin{equation*}
\Big|\alpha - \frac{a_1}{q_1} - \frac{a_2}{q_2} - ... -
\frac{a_n}{q_n} \Big| \ll_\epsilon \frac{1}{q N^{n - \epsilon}}
\end{equation*}
for some integers $a_1$, $a_2$, ..., $a_n$ and distinct prime
numbers $1 \leq q_1, q_2, ..., q_n \leq N$.
\end{conj}

Using Theorem \ref{thm1} with $n = 2m$ and $m$ very large in terms
of $\epsilon$, one can combine the first $m$ rational numbers to
$\frac{A_1}{Q_1}$ and the remaining to $\frac{A_2}{Q_2}$. Note that
$\kappa(n) = \frac{2n}{3} + O(1) = \frac{4m}{3} + O(1)$. Then one
can derive the following
\begin{cor} \label{cor1}
Let $\epsilon > 0$ and $N \geq 1$ be any real numbers. Suppose
$\alpha$ has a rational approximation $|\alpha - \frac{a}{q}| \leq
\frac{1}{q N^{4/3}}$ for some integers $a$, $1 \leq q \leq N^{4/3}$
and $(a,q) = 1$. Then
\begin{equation*}
\Big|\alpha - \frac{a_1}{q_1} - \frac{a_2}{q_2} \Big| \ll_\epsilon
\frac{1}{q N^{4/3 - \epsilon}}
\end{equation*}
for some integers $a_1$, $a_2$ and $1 \leq q_1, q_2 \leq N$.
\end{cor}

This is a bit weaker than Theorem \ref{thm0}. Meanwhile if we simply
combine the $n$ rational numbers to $\frac{A}{Q}$ where the number
of distinct prime factors of $Q$, $\omega(Q) = n$, we have
\begin{cor} \label{cor2}
For any real numbers $\epsilon > 0$ and $X \geq 1$, consider any
natural number $n \leq \sqrt{\frac{\epsilon \log X}{6 \log \log
X}}$. If $\alpha$ has a rational approximation $|\alpha -
\frac{a}{q}| \leq \frac{1}{q X}$ for some integers $a$, $1 \leq q
\leq X$ and $(a,q) = 1$, then there exists a rational number
$\frac{A}{Q}$ with $Q \leq X^{n / \kappa(n)}$ and $\omega(Q) = n$
such that
$$\Big| \alpha - \frac{A}{Q} \Big| \ll_\epsilon \frac{1}{q X^{1 -
\epsilon}}.$$ Note: When $n$ is large, $\frac{n}{\kappa(n)} \approx
\frac{3}{2}$.
\end{cor}
Proof of Corollary \ref{cor2}: Set $X = N^{\kappa(n)}$. We just need
to make sure the upper bound on $n$ in Theorem \ref{thm1} is
satisfied. Since $1/2 \leq \kappa(n) \leq 3n/4$, for $N$ or $X$
large enough,
$$n \leq \sqrt{\frac{\epsilon \log X}{6 \log \log X}} \Rightarrow n
\leq \sqrt{\frac{\epsilon n \log N}{6 \log \log N}} \Rightarrow
\sqrt{n} \leq \sqrt{\frac{\epsilon \log N}{6 \log \log N}}$$ which
gives the desired bound after squaring.

\bigskip

So, roughly speaking, given a rational approximation of $\alpha$
with denominator $\leq X$, we can always find a rational
approximation with denominator $\leq X^{3/2}$ having a prescribed
number of distinct prime divisors which approximates $\alpha$ nearly
as well. This leads to the following
\begin{conj} \label{conjecture}
For any real numbers $\epsilon > 0$ and $X \geq 1$, consider any
natural number $n \leq \epsilon \log X$. If $\alpha$ has a rational
approximation $|\alpha - \frac{a}{q}| \leq \frac{1}{q X}$ for some
integers $a$, $1 \leq q \leq X$ and $(a,q) = 1$, then there exists a
rational number $\frac{A}{Q}$ with $Q \leq X^{1 + \epsilon}$ and
$\omega(Q) = n$ such that
$$\Big| \alpha - \frac{A}{Q} \Big| \ll_\epsilon \frac{1}{q X^{1 -
\epsilon}}.$$
\end{conj}
{\bf Some Notations:} Throughout the paper, $\epsilon$ denotes a
small positive number. Both $f(x) = O(g(x))$ and $f(x) \ll g(x)$
mean that $|f(x)| \leq C g(x)$ for some constant $C > 0$. Moreover
$f(x) = O_\lambda(g(x))$ and $f(x) \ll_\lambda g(x)$ mean that the
implicit constant $C = C_\lambda$ may depend on the parameter
$\lambda$. Also $\omega(n)$ is the number of distinct prime factors
of $n$, and $e(x) = e^{2\pi i x}$. Finally $|\mathcal{S}|$ stands
for the cardinality of a set $\mathcal{S}$.
%-----------------------------------------------------------------------
\section{Some Preparations}
\begin{lem} \label{lem1}
Let $a$ and $q \geq 1$ be integers with $(a,q) = 1$. For any
positive integers $n$, $k$ and $L$, and any real number $N \geq 1$
such that
\begin{equation} \label{1cond}
L \leq N^n, \; q \leq L N^k, \; 2^{n+k+1} < N
\end{equation}
and $\mathcal{P}$ is the set of primes in $[N/2,N]$ not dividing
$q$, we have
\begin{equation} \label{1lem}
\sum_{l = 1}^{L} \Big| \sum_{q_1 \in \mathcal{P}} ... \sum_{q_n \in
\mathcal{P}} e\Bigl(l q_1 ... q_n \frac{a}{q} \Bigr) \Big| \ll 2^{n
+ k} n^n \max \Bigl(L N^{n/2 + k/2}, \frac{L N^n}{q^{1/2}} \Bigr).
\end{equation}
\end{lem}

Proof: Let $S$ be the right hand side of (\ref{1lem}). If $k \geq
n$, then trivially, $S \leq L N^n \leq L N^{n/2 + k/2}$. Now we can
assume $k < n$. Then
\begin{align} \label{s}
S =& \sum_{l=1}^{L} \sum_{q_1 \in \mathcal{P}} ... \sum_{q_k \in
\mathcal{P}} \Big| \sum_{q_{k+1} \in \mathcal{P}} ... \sum_{q_n \in
\mathcal{P}} e\Bigl(l q_1 ... q_n \frac{a}{q} \Bigr) \Big| \\
=& \sum_{(N/2)^k \leq r \leq L N^k} d_r \Big| \sum_{q_{k+1} \in
\mathcal{P}} ... \sum_{q_n \in \mathcal{P}} e\Bigl(r q_{k+1} ... q_n
\frac{a}{q} \Bigr) \Big| \nonumber
\end{align}
where
$$d_r = \mathop{\sum_{l=1}^{L} \sum_{q_1 \in \mathcal{P}} ... \sum_{q_k \in
\mathcal{P}}}_{l q_1 ... q_k = r} 1.$$ As $r \leq L N^k \leq N^{n +
k}$, $r$ is divisible by $m$ primes in $\mathcal{P}$ with $k \leq m
\leq n + k$ because $2^{n + k +1} < N$. Thus
$$d_r \leq \sum_{k \leq m \leq n+k} \binom{m}{k} k^k \leq k^k
\sum_{k \leq m \leq n+k} 2^m \leq 2^{n+k+1} k^k.$$ By Cauchy-Schwarz
inequality,
\begin{align} \label{s1}
S \leq& \Bigl(\sum_{(N/2)^k \leq r \leq LN^k} d_r^2 \Bigr)^{1/2}
\Bigl( \sum_{(N/2)^k \leq r \leq LN^k} \Big| \sum_{q_{k+1} \in
\mathcal{P}} ... \sum_{q_n \in \mathcal{P}} e\Bigl(r q_{k+1} ... q_n
\frac{a}{q} \Bigr) \Big|^2 \Bigr)^{1/2} \nonumber \\
\leq& 2^{n+k+1} k^k (LN^k)^{1/2} \Bigl( \sum_{r = 1}^{([LN^k/q]+1)q}
\Big| \sum_{q_{k+1} \in \mathcal{P}} ... \sum_{q_n \in \mathcal{P}}
e\Bigl(r q_{k+1} ... q_n \frac{a}{q} \Bigr) \Big|^2 \Bigr)^{1/2}
\nonumber \\
\leq& 2^{n+k+1} k^k \frac{LN^k}{q^{1/2}} \Bigl( \sum_{r = 1}^{q}
\sum_{q_{k+1}, ..., q_n \in \mathcal{P}} \sum_{q_{k+1}', ..., q_n'
\in \mathcal{P}} e\Bigl(r (q_{k+1} ... q_n - q_{k+1}' ... q_n')
\frac{a}{q} \Bigr) \Bigr)^{1/2}
\end{align}
where $[x]$ denotes the greatest integer less than or equal to $x$.
There are two cases. If $N^{n-k} \leq q$, then as an integer $Q =
q_{k+1} ... q_n$ has at most $(n-k)^{n-k}$ ways of representation,
the exponential sum in (\ref{s1}) is $\leq q N^{n-k} (n - k)^{n -
k}$ using
\begin{equation} \label{ortho}
\sum_{r=1}^{q} e\Bigl(\frac{a m r}{q}\Bigr) = \Bigl\{
\begin{tabular}{ll}$q$, & if $q$ divides $m$, \\
$0$, & otherwise.
\end{tabular}
\end{equation}
Putting this into (\ref{s1}), we have $S \ll 2^{n+k} n^n L N^{n/2 +
k/2}$.

If $q < N^{n-k}$, we divide the interval $[1, N^{n-k}]$ into
$N^{n-k}/q + O(1)$ intervals $I_i$ of length at most $q$. Then the
exponential sum in (\ref{s1}) is
\begin{align*}
=& \sum_{i} \sum_{j} \sum_{r = 1}^{q} \sum_{q_{k+1} ... q_n \in I_i}
\sum_{q_{k+1}' ... q_n' \in I_j} e\Bigl(r (q_{k+1} ... q_n -
q_{k+1}' ... q_n') \frac{a}{q} \Bigr) \\
\ll& \Bigl(\frac{N^{n-k}}{q}\Bigr)^2 q^2 (n-k)^{(n-k)}
\end{align*}
by (\ref{ortho}). Putting this into (\ref{s1}), $S \ll 2^{n+k} n^n L
N^n / q^{1/2}$. Combining the two cases, we have the lemma.
%-----------------------------------------------------------------------
\begin{lem} \label{lem2}
Assume the same conditions for $a$, $q$, $n$, $k$, $L$, $N$ and the
same $\mathcal{P}$ as in Lemma \ref{lem1}. Suppose further that we
have a positive integer $m \leq n$ and $r_1 + ... + r_{m} = n$ with
positive integers $r_1$, ... $r_{m}$. Then
\begin{equation} \label{2lem}
S = \sum_{l = 1}^{L} \Big| \sum_{q_1 \in \mathcal{P}} ...
\sum_{q_{m} \in \mathcal{P}} e\Bigl(l q_1^{r_1} ... q_{m}^{r_{m}}
\frac{a}{q} \Bigr) \Big| \ll 2^{n + k} n^n \max \Bigl(L N^{n/2 +
k/2}, \frac{L N^n}{q^{1/2}} \Bigr).
\end{equation}
\end{lem}

Proof: If $m \leq n/2 + k/2$, then we have trivially that $S \leq L
N^{m} \leq L N^{n/2 + k/2}$. If $m > n/2 + k/2$, let $M$ be the
number of $1$'s among $r_1$, ..., $r_{m}$. Then $M + 2 (m - M) \leq
r_1 + ... + r_{m} = n$ which gives $M \geq 2m - n > 2(n/2 + k/2) - n
= k$. So there are $1 \leq i_1 < ... < i_k \leq m$ such that
$r_{i_1} = ... = r_{i_k} = 1$ and one can split $q_{i_1}$, ...
$q_{i_k}$ out similar to (\ref{s}) in Lemma \ref{lem1}. Hence one
can imitate the proof in Lemma \ref{lem1} and get the same upper
bound.

\bigskip

We also need the following M\"{o}bius inversion formula (see
[\ref{L}, \S25] for example).
\begin{lem} \label{lem3}
$$\mathop{\sum_{q_1} ... \sum_{q_n}}_{\hbox{distinct}} f(q_1, ...,
q_n) = \sum_{\mathcal{S}} \mu(\mathcal{S}) \mathop{\sum_{q_1} ...
\sum_{q_n}}_{q_i \hbox{ with restriction in } \mathcal{S} } f(q_1,
..., q_n)$$ where $\mathcal{S}$ is over all possible partition of
$\{1, 2, ..., n\}$, say $\mathcal{S} = \{P_1, ..., P_m\}$ such that
$\bigcup_{i = 1}^{m} P_i = \{1, ..., n \}$, and $\mu(\mathcal{S}) =
\prod_{j = 1}^{m} (-1)^{|P_j| - 1} (|P_j| - 1)!$. The $q_i$ with
restriction in $\mathcal{S}$ means the requirement $q_i = q_j$ for
$i, j \in P_k$, $1 \leq k \leq m$.
\end{lem}

%-------------------------------------------------------------------------
\section{Proof of Theorem \ref{thm1}}

First, we recall a variant of Erd\"{o}s - Tur\'{a}n inequality on
uniform distribution (see, for example, R.C. Baker [\ref{B}, Theorem
2.2]).
\begin{lem} \label{baker}
Let $L$ and $J$ be positive integers. Suppose that $||x_j|| \geq
\frac{1}{L}$ for $j= 1, 2,... , J$. Then
$$\sum_{l = 1}^{L} \Big| \sum_{j = 1}^{J} e(l x_j) \Big| >
\frac{J}{6}.$$ Here $||x|| = \min_{n \in \mathbb{Z}} |x - n|$, the
distance from $x$ to the nearest integer.
\end{lem}

Proof of Theorem \ref{thm1}: Let $\epsilon > 0$, $N \geq 1$ and $n/2
\leq \phi \leq n$. Let $\mathcal{P}$ be the set of prime numbers in
the interval $[N/2, N]$ not dividing $q$. Without loss of
generality, we can suppose that $\alpha$ has a rational
approximation $|\alpha - \frac{a}{q}| \leq \frac{1}{q N^{\phi}}$ for
some integers $a$, $N^{n - \phi} \leq q \leq N^{\phi}$ and $(a,q) =
1$. The case where $q \leq N^{n - \phi}$ is covered by the trivial
bound as one can simply choose some prime numbers $N/2 \leq q_1 <
q_2 < ... < q_n \leq N$, and there exists integer $b$ such that
$$\Big|\alpha - \frac{b}{q_1 q_2 ... q_n}\Big| \leq \frac{n+1}{q_1 q_2 ... q_n}
\leq \frac{(n+1) 2^n}{N^n} \ll_\epsilon \frac{1}{N^{n-\epsilon}}
\leq \frac{1}{q N^{\phi - \epsilon}}$$ for $N$ large in terms of
$\epsilon$. We are trying to find integer $b$ and distinct prime
numbers $q_1, q_2, ..., q_n \in \mathcal{P}$ such that
\begin{equation} \label{S0}
\Big| \frac{a}{q} - \frac{b}{q_1 q_2 ... q_n} \Big| < \frac{1}{q
N^{\phi - \epsilon/2}} \; \hbox{ or } \; \Big\Vert \frac{q_1 q_2 ...
q_n a}{q} \Big\Vert < \frac{1}{q N^{\phi - n - \epsilon/2}}.
\end{equation}
In view of the Lemma \ref{baker}, it suffices to show
\begin{equation} \label{S}
S := \sum_{l = 1}^{L} \Big| \mathop{\sum_{q_1, q_2, ..., q_n \in
\mathcal{P}}}_{\hbox{distinct}} e \Bigl(\frac{l q_1 q_2 ... q_n
a}{q} \Bigr) \Big| \leq \frac{1}{6} (|\mathcal{P}|^n - n^2
|\mathcal{P}|^{n-1})
\end{equation}
with $L := [q N^{\phi - n}] + 1$ because there are at least
$|\mathcal{P}|^n - \binom{n}{2} |\mathcal{P}|^{n-1}$ numbers of the
form $q_1 q_2 ... q_n$ with $q_i$ distinct. Now we apply Lemma
\ref{lem3} to (\ref{S}), and bound it by using Lemma \ref{lem2},
$|\mu(\mathcal{S})| \leq n!$, an upper bound for the number of
partitions and Stirling's formula,
\begin{align*}
S =& \sum_{l = 1}^{L} \Big| \sum_{\mathcal{S}} \mu(\mathcal{S})
\mathop{\sum_{q_1} ... \sum_{q_n}}_{q_i \hbox{ with restriction in }
\mathcal{S} } e \Bigl(\frac{l q_1 q_2 ... q_n a}{q} \Bigr) \Big| \\
\leq& \sum_{\mathcal{S}} |\mu(\mathcal{S})| \sum_{l = 1}^{L} \Big|
\mathop{\sum_{q_1} ... \sum_{q_n}}_{q_i \hbox{ with restriction in }
\mathcal{S} } e \Bigl(\frac{l q_1 q_2 ... q_n a}{q} \Bigr) \Big| \\
\ll& e^{\pi \sqrt{2/3} \sqrt{n}} \sqrt{2 \pi n} \Bigl(\frac{n}{e}
\Bigr)^n 2^{n + k} n^n \max \Bigl(L N^{n/2 + k/2}, \frac{L
N^n}{q^{1/2}} \Bigr) \\
\ll& (2n)^{2n} \max \Bigl(L N^{n/2 + k/2}, \frac{L N^n}{q^{1/2}}
\Bigr)
\end{align*}
for any integer $k \geq n - \phi$. Note: One can check that
condition (\ref{1cond}) is satisfied for our choice of $L$ and $k$.
As $n^2 \leq (\log N)^2$, the bound in (\ref{S}) holds if
\begin{equation} \label{S1}
(2n)^{2n} \max \Bigl(L N^{n/2 + k/2}, \frac{L N^n}{q^{1/2}} \Bigr)
\leq \frac{N^n}{(3 \log N)^{n+1}}
\end{equation}
for $N$ sufficiently large, as $|\mathcal{P}| \geq \frac{N}{3 \log
N}$. In other words, we need
\begin{align}
12^n n^{2n} (\log N)^{n+1} q N^{\phi - n} N^{n/2 + k/2} \ll N^n,
\label{c1} \\
12^n n^{2n} (\log N)^{n+1} q^{1/2} N^{\phi - n} N^{n} \ll N^n.
\label{c2}
\end{align}
After some algebra and using $q \leq N^{\phi}$, (\ref{c1}) and
(\ref{c2}) are true when
\begin{align*}
\phi < \frac{3n}{4} - \frac{k}{4} - \frac{n \log 4n}{\log N} -
\frac{(n+1) \log \log N}{2 \log N}, \\
\phi < \frac{2n}{3} - \frac{4n \log{4 n}}{3 \log N} -
\frac{2(n+1) \log \log N}{3 \log N}.
\end{align*}
Thus, (\ref{S}) is true when
$$\phi \geq n - k, \; \phi < \frac{3n}{4} - \frac{k}{4} - 2\epsilon'
\hbox{ and } \phi < \frac{2n}{3} - 3\epsilon '$$ as $n \leq
\frac{\epsilon' \log N}{\log \log N}$. Set $k = [n/3] + 1$, one can
check that
$$n - [n/3] - 1 < \frac{3n}{4} - \frac{[n/3] + 1}{4} - 3\epsilon' <
\frac{2n}{3} - 3\epsilon'$$ for $\epsilon'$ sufficiently small.
Therefore, we can choose $\phi = \frac{3n}{4} - \frac{[n/3] + 1}{4}
- 3\epsilon'$ and this gives the theorem by setting $\epsilon' =
\epsilon/6$.
%-----------------------------------------------------------------------
\section{Proof of Theorem \ref{thm2}}

First we quote a lemma (see Chapter 1, Lemma 8a of [\ref{V}]).
\begin{lem} \label{vin}
Let $a$ and $q \neq 0$ be relatively prime integers and suppose $N
\geq 1$. Then
$$\sum_{r = 1}^{q} \min \Bigl(N, \frac{1}{||\frac{a r}{q}||} \Bigr)
\ll N + q \log q.$$
\end{lem}

Proof of Theorem \ref{thm2}: It starts the same way as the proof of
Theorem \ref{thm1} with $n = 3$. Suppose we consider only $3/2 \leq
\phi \leq 2$. Then $N \leq q \leq N^\phi \leq N^2$. Now, instead of
having all the denominators $q_1$, $q_2$, $q_3$ in $\mathcal{P}$, we
allow $q_3$ to be simply from the interval $[N/2 , N]$ but different
from $q_1$ and $q_2$ as used in [\ref{C}]. Borrowing from (\ref{S}),
it suffices to prove
\begin{equation} \label{S2}
S := \sum_{l = 1}^{L} \Big| \mathop{\sum_{q_1 \in \mathcal{P}}
\sum_{q_2 \in \mathcal{P}} \sum_{N/2 \leq q_3 \leq N}
}_{\hbox{distinct}} e \Bigl(\frac{l q_1 q_2 q_3 a}{q} \Bigr) \Big|
\leq \frac{1}{6} (|\mathcal{P}|^3 - n^2 |\mathcal{P}|^{2})
\end{equation}
with $L := [q N^{\phi - 3}] + 1$ because there are at least
$|\mathcal{P}|^3 - \binom{n}{3} |\mathcal{P}|^{2}$ numbers of the
form $q_1 q_2 q_3$ with $q_i$ distinct primes. Splitting out $q_1$
and $q_2$, we have
\begin{equation*}
\begin{split}
S \leq& \sum_{l = 1}^{L} \sum_{q_1 \neq q_2 \in \mathcal{P}} \Big|
\mathop{\sum_{N/2 \leq q_3 \leq N}}_{q_3 \neq q_1, q_2} e
\Bigl(\frac{l q_1 q_2 q_3 a}{q} \Bigr) \Big| \\
\leq& \binom{5}{2} \sum_{r = 1}^{L N^2} \Big| \sum_{N/2 \leq q_3
\leq N} e \Bigl(\frac{r q_3 a}{q} \Bigr) \Big| + 2 LN^2
\end{split}
\end{equation*}
as $LN^2 \leq 2 q N^{\phi - 3} N^2 \leq 2 N^5$ and $r$ can be
divisible by at most $5$ primes in $\mathcal{P}$. The summation
above is
\begin{equation*}
\begin{split}
\ll& \sum_{r = 1}^{q([L N^2/q] + 1)} \Big| \sum_{N/2 \leq q_3 \leq
N} e \Bigl(\frac{r q_3 a}{q} \Bigr) \Big| \leq \Bigl( [\frac{L
N^2}{q}] + 1 \Bigr) \sum_{r = 1}^{q} \Big| \sum_{N/2 \leq q_3
\leq N} e \Bigl(\frac{r q_3 a}{q} \Bigr) \Big| \\
\ll& \frac{L N^2}{q} \sum_{r = 1}^{q} \min \Bigl( N,
\frac{1}{||\frac{a r}{q}||} \Bigr) \ll \frac{L N^2}{q} (N + q \log
q) \ll LN^2 \log q
\end{split}
\end{equation*}
by Lemma \ref{vin} as $N \leq q$. Thus $S \ll L N^2 \log q$ and
(\ref{S2}) is true if $L N^{2 + \epsilon'} \ll N^3$ for any
$\epsilon' > 0$. In other words, we want $q N^{\phi - 3} \ll N^{1 -
\epsilon'}$. So $q \ll N^{4 - \phi - \epsilon'}$ and on the other
hand we require $q \leq N^\phi$. The optimum choice for $\phi = 2 -
\epsilon'$ and this gives the theorem after putting $\phi$ back into
(\ref{S0}).

%-----------------------------------------------------------------------

\bigskip
Department of Mathematical Sciences \\
University of Memphis \\
Memphis, TN 38152 \\
U.S.A. \\
tszchan@memphis.edu
\end{document}